\documentclass[notitlepage]{amsart}

\newtheorem{pro}{Proposition}[section]
\newtheorem{teo}[pro]{Theorem}
\newtheorem{defi}[pro]{Definition}
\newtheorem{lem}[pro]{Lemma}
\newtheorem{cor}[pro]{Corollary}
\newtheorem{rk}[pro]{Remark}

\newcommand{\Ext}{\mathrm{Ext}}
\newcommand{\Hom}{\mathrm{Hom}}

\newcommand{\A}{\mathcal{A}}
\newcommand{\B}{\mathcal{B}}
\newcommand{\I}{\mathcal{I}}

\newcommand{\Q}{\mathcal{P}}
\newcommand{\X}{\mathcal{X}}
\newcommand{\Y}{\mathcal{Y}}

\newcommand{\pd}{\mathrm{pd}}
\newcommand{\id}{\mathrm{id}}
\newcommand{\resdim}{\mathrm{resdim}}
\newcommand{\Findim}{\mathrm{Findim}\,}
\newcommand{\findim}{\mathrm{findim}\,}
\newcommand{\coresdim}{\mathrm{coresdim}}

\newcommand{\Add}{\mathrm{Add}}
\newcommand{\Ker}{\mathrm{Ker}}
\newcommand{\Coker}{\mathrm{Coker}}
\newenvironment{dem}{\noindent\bf Proof. \rm }{$\ \Box$}
\usepackage{latexsym,amssymb,amscd} 
\usepackage{amsmath}

\begin{document}

\title{Homological Dimensions in Cotorsion Pairs}
\author{Lidia Angeleri H$\ddot{\mathrm{u}}$gel and Octavio Mendoza}

\begin{abstract} 
Two classes $\mathcal A$ and $\mathcal B$  of modules over a ring $R$ are said to form a cotorsion pair $(\mathcal A, \mathcal B)$ if $\mathcal A={\rm Ker\,Ext}^1_R(-,\mathcal B)$ and 
$\mathcal B={\rm Ker\,Ext}^1_R(\mathcal A,-)$.
We investigate  relative homological dimensions in cotorsion pairs. This can be applied to study the big and the little finitistic dimension of $R$. We show that $\Findim R<\infty$ if and only if the following dimensions are finite for some  cotorsion pair $(\mathcal A, \mathcal B)$ in $\mathrm{Mod}\,R$:  the relative projective dimension of $\A$ with respect to itself,  and the $\mathcal A$-resolution dimension of the category $\mathcal P$ of all $R$-modules of finite projective dimension. Moreover, we obtain an analogous result for $\findim R$, and we characterize when $\Findim R=\findim R.$
\end{abstract}  
\maketitle
\section*{Introduction.}
 The   study  of homological dimensions which are obtained by replacing the projective or injective modules by certain subcategories  was initiated by Auslander and Buchweitz in their seminal paper \cite{AuB}, which was the starting point for what is now called relative homological algebra.
 \
 
  Of course, the existence of approximations is the prerequisite for computing relative dimensions. In recent years, a powerful machinery for producing approximations was developed by using the notion of a cotorsion pair, see \cite{ET, EJ,  GT}. So it is not surprising that
 cotorsion pairs provide a good setting for investigating relative   homological dimensions.
 \

The main purpose of this note is to use cotorsion pairs in order to obtain estimates for the finitistic dimensions.
 Recall that the   (right)  \emph{big finitistic dimension} of $R$ is defined as $\Findim(R_R)=\pd\,\Q,$ the supremum of the projective dimensions attained on the subcategory $\Q$ of all modules of finite projective dimension. By restricting to the subcategory $\Q^{<\infty}$ of all modules in $\Q$ that admit a projective resolution consisting of finitely generated projective modules, one obtains    $\findim(R_R)=\pd\,\Q^{<\infty},$    the (right)  \emph{little finitistic dimension} of $R$.
The Finitistic Dimension Conjecture is one of the main  open problems in the representation theory of algebras. It asserts that $\findim(R_R)<\infty$ for any artin algebra $R$.
\

\

Let now $(\A,\B)$ be a complete hereditary cotorsion pair (see the definitions at the end of Section 1), and assume that the relative projective dimension
$$\pd_\A\,(\A)=\mathrm{min}\,\{n\geq 0\,:\,\Ext_R^j(,-)|_{\A}=0
  \text{ for any } j>n\}$$
  is finite. Denote by $\resdim_\A$ the resolution dimension with respect to $\A$. We prove the following inequality.
  \
  
  \
  
  {\bf Theorem \ref{Thmprincipal}}
  $$ \resdim_\A\,(\Q^{<\infty})\leq\findim(R_R)\leq\pd_\A\,(\A)+\resdim_\A\,(\Q^{<\infty}).$$
\

  As a consequence, we obtain a criterion for validity of the Finitistic Dimension Conjecture.
  \
  
  \
  
  {\bf Corollary \ref{equivFindim}}
  {\it The following conditions are equivalent.
 \begin{itemize}
  \item[(a)] $\findim(R_R)<\infty.$
  \item[(b)] There is a hereditary complete cotorsion pair $(\A, \B)$ in $\mathrm{Mod}\,R$ such that $\pd_\A\,(\A)<\infty$ and $ \resdim_\A\,(\Q^{<\infty})<\infty.$
  \item[(c)] Every hereditary complete cotorsion pair $(\A, \B)$ in $\mathrm{Mod}\,R$ with $\pd_\A\,(\A)<\infty$ satisfies $ \resdim_\A\,(\Q^{<\infty})<\infty.$
 \end{itemize} }
 
 \medskip
 
 An analogous result characterizes finiteness of $\Findim(R_R)$. Moreover, we characterize equality of the big and the little finitistic dimension in terms of the kernel $\omega=\A\cap\B$ of a certain cotorsion pair $(\A, \B)$, see Theorem \ref{Findim-findim}.
This  generalizes a result from \cite{AnT1}.
\

\

 Finally, we prove that the finitistic dimension is bounded by the homological dimensions of tilting modules. 
 \
 
 \
 
{\bf Corollary \ref{pd+id}} Let $T$ be a tilting module in $\mathrm{Mod}\,R.$ If $R$ is right noetherian then $$\Findim(R_R)\leq\pd\,T+\id\,T.$$

\

As an application, we improve a result of  V. Mazorchuk and S. Ovsienko  \cite{MO} on properly stratified algebras having a simple preserving duality. More precisely,  let $R$ be such an algebra, and  assume  that every classical tilting right $R$-module is also cotilting. Denote by $T$   the   characteristic tilting module associated to   $R.$   Then it follows from Corollary \ref{pd+id} that $\Findim(R_R)\le 2\,\pd\,T,$     see Remark \ref{applications}. Combining this with  \cite[Theorem 1]{MO} we conclude that
$\Findim(R_R)=\findim(R_R)=2\,\pd\,T$.

\bigskip

This paper is  organized as follows. Section 1 is devoted to  some preliminaries. In Section 2, we use   results from \cite{MS} in order to obtain  formulae that relate relative projective dimensions and resolution dimensions in cotorsion pairs.  In Section 3, we apply these formulae to cotorsion pairs arising from tilting theory. Hereby we refine some results from \cite{AnC} and \cite{AnHT}. The finitistic dimensions are the topic of Section 4.  
 
 \bigskip

\section{Preliminaries}
\

We start this section by collecting all the background material that
will be necessary in the sequel. First, we
introduce some general notation. Next, we recall the notion
of relative projective dimension and resolution dimension of a given
class of modules. Finally, we also recall  definition and basic properties of cotorsion pairs.
\

Throughout the paper, $R$ will be an associative ring with unit, and $\mathrm{Mod}\,R$ the category
of all right $R$-modules. The subcategory of all modules
possessing a projective resolution consisting of finitely generated modules is
denoted by $\mathrm{mod}\,R.$ By a subcategory of 
$\mathrm{Mod}\,R$ we always mean a full subcategory. 
\

\

We denote by $\mathrm{pd}\,X$ the  {\bf{projective dimension}} of $X.$ Similarly, $\mathrm{id}\,X$ is the {\bf{injective dimension}} of $X.$ 
For any non-negative integer $n,$ we set 
$$\Q_n:=\{X\in\mathrm{Mod}\,R\;:\;\pd X\leq n\}$$ and $\Q^{<\infty}_n:=\Q_n\cap\mathrm{mod}\,R;$ moreover, $$\Q:=\{X\in\mathrm{Mod}\,R\;:\;\pd X<\infty\}$$ and $\Q^{<\infty}:=\Q\cap\mathrm{mod}\,R.$ The classes $\I,$ $\I_n,$ $\I^{<\infty}$ and $\I^{<\infty}_n$ are defined dually.  In particular, $\Q_0$ consists of  the projective $R$-modules and $\I_0$ of the injectives. 
\

\

Let now  $\X$ be a subcategory of $\mathrm{Mod}\,R$. We denote by $\Add\,(\X)$   the class of all  $R$-modules isomorphic to direct summands of direct sums of elements in $\mathcal{X}.$ Moreover, for each positive integer $i,$ we denote $$\X^{\perp_i}:=\{M\in\mathrm{Mod}\,R\;:\;\Ext^i_R(-,M)|_{\X}=0\}\quad {\rm and}\quad  \X^\perp:=\cap_{i>0}\,\X^{\perp_i}.$$ Dually, we have the classes ${}^{\perp_i}\X$ and ${}^{\perp}\X.$
\ 

\

A subcategory $\X\subseteq\mathrm{Mod}\,R$ (respectively, $\X\subseteq\mathrm{mod}\,R$)  is said to be {\bf{resolving}} if it is closed under extensions and kernels of surjections, and it contains all (respectively, all finitely generated) projective modules.  If the dual properties hold true, then $\X$ is a {\bf{coresolving}} subcategory. 
For example, $\Q$  and ${}^{\perp}\X$ are resolving subcategories of $\mathrm{Mod}\,R$, while $\I$ and $\X^\perp$ are coresolving subcategories of $\mathrm{Mod}\,R$.

\

\bigskip

{\sc Relative homological dimensions.} Given a class of $R$-modules $\X$ and an $R$-module $M,$ the {\bf{relative projective dimension}} of $M$ with respect to $\X$ is defined as $$\pd_{\X}\,(M):=\mathrm{min}\,\{n\geq 0\,:\,\Ext_R^j(M,-)|_{\X}=0  \text{ for any } j>n\}.$$ Dually, we denote by
  $\mathrm{id}_{\X}\,(M)$ the  {\bf{relative injective dimension}} of
  $M$ with respect to $\X.$ Furthermore, for any class $\Y\subseteq\mathrm{Mod}\,R,$ we set $$\pd_\X\,(\Y):=\mathrm{sup}\,\{\pd_\X\,Y\;:\; Y\in\Y\}\text{ and }\id_\X\,(\Y):=\mathrm{sup}\,\{\id_\X\,Y\;:\; Y\in\Y\}.$$ 
  If $\X=\mathrm{Mod}\,R,$ we just write $\pd\,(\Y)$ and   $\id\,(\Y)$.
 \bigskip
 
   The following basic properties are straightforward.

\begin{lem}\label{basic} Let $\X$ and $\Y$  be  classes in $\mathrm{Mod}\,R.$ Then $$\pd_\X\,(\Y)=\id_\Y\,(\X).$$ Furthermore, for any exact sequence  $0\to M'\to M\to M''\to 0$ in $\mathrm{Mod}\,R$ we have  
 \begin{itemize}
  \item[(a)] $\pd_\X\,(M)\leq\mathrm{max}\,\{\pd_\X\,(M'),\,\pd_\X\,(M'')\};$
\vspace{.2cm}
  \item[(b)] $\pd_\X\,(M')\leq\mathrm{max}\,\{\pd_\X\,(M),\,\pd_\X\,(M'')-1\};$
\vspace{.2cm}
  \item[(c)] $\pd_\X\,(M'')\leq\mathrm{max}\,\{\pd_\X\,(M),\,\pd_\X\,(M')+1\}.$
 \end{itemize} 
\end{lem}

\bigskip

{\sc Resolution and coresolution dimension.}
Let $M\in\mathrm{Mod}\,R$ and $\X$ be a class of $R$-modules. The $\X$-{\bf{coresolution dimension}} $\coresdim_\X\,(M)$ of $M$ is the minimal non-negative integer $n$ such that there is an exact sequence $$0\to M\to X_0\to X_1\to\cdots\to X_n\to 0$$ with $\X_i\in\X$ for $0\leq i\leq n.$ If such $n$ does not exist, we set $\coresdim_\X\,(M):=\infty.$ Also, we denote by $\mathcal{X}^{\vee}$ the  class of $R$-modules having finite $\mathcal{X}$-coresolution. 

Dually, we have the $\X$-{\bf{resolution dimension}} $\resdim_\X\,(M)$ of $M,$ and the class $\mathcal{X}^{\wedge}$ of $R$-modules having finite $\mathcal{X}$-resolution. 

Given a class $\Y\subseteq\mathrm{Mod}\,R,$ we set  
$$\coresdim_\X\,(\Y):=\mathrm{sup}\,\{\coresdim_\X\,(Y)\;
:\; Y\in\Y\},$$ and  $\resdim_\X\,(\Y)$ is defined dually.
\

\bigskip

{\sc Approximations.}  Let  $\X$ be a class of $R$-modules. A morphism $f:X\rightarrow M$ is said to be an $\X$-{\bf{precover}} if $X\in\X$ and $\mathrm{Hom}_R(Z,f):\mathrm{Hom}_R(Z,X)\rightarrow\mathrm{Hom}_R(Z,M)$
is surjective for any $Z\in\X.$ Furthermore, an $\X$-precover $f:X\rightarrow M$ is {\bf{special}} if $\Ker\,(f)=0$ and $\Coker\,(f)\in{}^{\perp_1}\X.$ 
We will  freely use also the dual notion of (special) $\X$-{\bf{preenvelope}}.
\

\

Finally, we recall  the notion of cotorsion pair which was introduced by L. Salce in \cite{S}. It is the analog of a torsion pair where the functor $\Hom_R(-,-)$ is replaced by $\Ext^1_R(-,-).$

\begin{defi}\cite{S} Let $\A$ and $\B$ be classes in $\mathrm{Mod}\,R.$ The pair $(\A,\B)$ is said to be a {\bf{cotorsion pair}} if $\A={}^{\perp_1}\B$ and $\A^{\perp_1}=\B.$ The class $\A\cap\B$ is called the {\bf{kernel}} of the cotorsion pair $(\A,\B)$.
\end{defi}

\begin{lem} \cite[Corollary 2.4]{S} For a cotorsion pair $(\A,\B)$ in $\mathrm{Mod}\,R,$ the following conditions are equivalent.
 \begin{itemize}
  \item[(a)] Every $R$-module has a special $\A$-precover.
\vspace{.2cm}
  \item[(b)] Every $R$-module has a special $\B$-preenvelope.
 \end{itemize}
In this case, the cotorsion pair $(\A,\B)$ is called {\bf{complete}}.
\end{lem}

\begin{lem} For a cotorsion pair $(\A,\B)$ in $\mathrm{Mod}\,R,$ the following conditions are equivalent.
 \begin{itemize}
  \item[(a)] $\A$ is resolving.
\vspace{.2cm}
  \item[(b)] $\B$ is coresolving.
\vspace{.2cm}
  \item[(c)] $\id_\A\,(\B)=0.$
 \end{itemize}
In this case, the cotorsion pair $(\A,\B)$ is called {\bf{hereditary}}.
\end{lem}
\bigskip

\section{General Results}

In this section, we obtain general results for complete cotorsion pairs. They will later be applied, on one hand, to tilting theory, and on the other hand, to the  big  and the  little  finitistic dimension of the ring $R.$ 

Let us start by collecting some preliminary results. They are stated in \cite{MS} under stronger assumptions, but the  proofs there also  work  in the present context. 

\begin{pro} \cite[Theorem 2.1 and Lemma 3.3]{MS} \label{Prop1} Let $\X$ and $\Y$ be classes in $\mathrm{Mod}\,R.$  The following statements hold true.  
 \begin{itemize}
\vspace{.2cm}
  \item[(a)]  $\id_\X\,(L)\leq$ $\id_\X\,(\Y) + \coresdim_\Y\,(L)$  for every $R$-module $L.$ 
\vspace{.2cm}
  \item[(b)] Assume that $\Y=\X^{\perp_1}$, or that $\Y$ is a subcategory of $\X$ which is closed under direct summands. Suppose further that $\id_\X\,(\Y)=0$. Then we have  \linebreak
$\id_\X\,(L)= \coresdim_\Y\,(L)$ for any $L\in\Y^\vee.$
\vspace{.2cm}
  \item[(c)] Let $\Y=\X^{\perp_1}$. 
 Then $\coresdim_{\Y}\,(M)\leq\id_\X(M)$ for any $M\in\mathrm{Mod}\,R.$ 
 
In particular, $\coresdim_{\Y}\,(\mathrm{Mod}\,R)\leq\pd\,\X.$
 \end{itemize}
\end{pro}

Of course, these results also have dual versions which we will freely use in the sequel.
We now turn to complete cotorsion pairs.

\begin{teo}\label{IneqDim} Let $(\A,\B)$ be a complete cotorsion pair in $\mathrm{Mod}\,R.$
 \begin{itemize}
 \item[(a)]  For any $R$-module $M,$ we have  $$\id\,M=\mathrm{max}\,\{\id_\A\,(M),\id_\B\,(M)\}$$
$$\pd\,M=\mathrm{max}\,\{\pd_\A\,(M),\pd_\B\,(M)\}.$$
  \item[(b)]  $\coresdim_\B\,(\mathrm{Mod}\,R)\leq\pd\,\A\leq\id_\A\,(\B)+ \coresdim_\B\,(\A)+1.$
\vspace{.2cm}
  \item[(c)] $\resdim_\A\,(\mathrm{Mod}\,R)\leq\id\,\B\leq\pd_\B\,(\A)+ \resdim_\A\,(\B)+1.$
  \end{itemize}
\end{teo}
\begin{dem} 
(a) We only prove the first equality, the second one follows similarly. It is clear that $\mathrm{max}\,\{\id_\A\,(M),\id_\B\,(M)\}\leq\id\,M.$\\
Take an $R$-module $N.$ Since $(\A,\B)$ is complete, we have an exact sequence $$0\to N\to B\to A\to 0$$ with $B\in\B$ and $A\in\A.$ By \ref{basic} (b), we infer that
$$\pd_{\{M\}}(N)\leq\mathrm{max}\,\{\pd_{\{M\}}(\B),\pd_{\{M\}}(\A)\}$$ for any $N\in\mathrm{Mod}\,R,$ which proves the result since $\pd_{\{M\}}(\mathrm{Mod}\,R)=\id\,M.$
\\
(b) The first inequality   follows from \ref{Prop1} (c) since $\B=\A^{\perp_1}.$ To prove the second one, let $M$ be an $R$-module. Since $(\A,\B)$ is complete, there is an exact sequence $0\to M\to B\to A\to 0$ with $A\in\A$ and $B\in\B.$ Thus, by the dual of \ref{basic} (c), we get that $$\id_\A(M)\leq\mathrm{max}\,\{\id_\A\,(\B),\,\id_\A\,(\A)+1\}\;\text{ for any }\;M\in\mathrm{Mod}\,R.$$ Furthermore, we know from  \ref{Prop1} (a) that   $\id_\A\,(\A)\leq$ $\id_\A\,(\B) + \coresdim_\B\,(\A)$, so the proof is complete. 
\\
(c) It is shown similarly.
\end{dem}

\begin{cor}\label{equa-ineq} Let $(\A,\B)$ be a hereditary complete cotorsion pair in $\mathrm{Mod}\,R.$ Then
 \begin{itemize}
  \item[(a)] $\pd_\A\,(\A)=\pd\,\A$ and $\id_\B\,(\B)=\id\,\B.$
\vspace{.2cm}
  \item[(b)] $\coresdim_\B\,(\mathrm{Mod}\,R)\leq\pd\,\A\leq \coresdim_\B\,(\A)+1.$
\vspace{.2cm}
  \item[(c)] $\resdim_\A\,(\mathrm{Mod}\,R)\leq\id\,\B\leq \resdim_\A\,(\B)+1.$
  \end{itemize}
\end{cor}
\begin{dem} It follows easily from Theorem \ref{IneqDim} since $\pd_\B\,(\A)=\id_\A\,(\B)=0.$
\end{dem}

\

The following result  will be useful in investigating  cotorsion pairs where $\pd\,\A$ is finite. It relies on work of M. Auslander and R. O. Buchweitz.
\

\begin{lem}\label{equal} Let $(\A,\B)$ be a hereditary complete cotorsion pair in $\mathrm{Mod}\,R$ with kernel $\omega.$ The following statements hold true.
 \begin{itemize}
  \item[(a)] $\A\cap\omega^\vee=\{X\in\A\,:\,\id_\A\,(X)<\infty\}.$
\vspace{.2cm}
 \item[(b)] $\pd_\B\,(M)=\resdim_\omega\,(M)$ for any $M\in\omega^\wedge.$
\vspace{.2cm}
  \item[(c)] $\pd_\omega\,(M)=\resdim_\A\,(M)$ for any $M\in\A^\wedge.$
 \end{itemize}
\end{lem}
\begin{dem} First of all, we observe  that $\omega$ is an injective cogenerator for $\A$ in the sense of Auslander-Buchweitz  \cite[Pag. 17]{AuB}, since $(\A,\B)$ is hereditary and complete. Thus statements (a) and (c)
 follow immediately from \cite[Lemma 4.3 and Proposition 2.1]{AuB}.
As for (b), it a consequence of the dual result of  \ref{Prop1} (b), because $\omega$ is a subcategory of $\B$ closed under direct summands with $\pd_\B\,(\omega)=0.$
\end{dem}

\

Now we are ready to state one of our key results.
\

\begin{teo}\label{Thmequal} Let $(\A,\B)$ be a hereditary complete cotorsion pair in $\mathrm{Mod}\,R$ with kernel $\omega.$ The following statements hold true.
 \begin{itemize}
  \item[(a)] $\pd\,\A=\pd_\A\,(\A)=\coresdim_\B\,(\A)=\coresdim_\omega\,(\A)=\coresdim_\B\,(\mathrm{Mod}\,R).$
\vspace{.2cm}
  \item[(b)] $\pd_\A\,(\A)<\infty$ if and only if $\A\subseteq\omega^\vee$ and $\pd\,\omega<\infty.$
In this case, we have $$\pd_\A\,(\A)=\pd\,\omega=\id_\omega\,(\A).$$
 \end{itemize}
\end{teo}
\begin{dem} 
First of all, recall that $\pd_\A\,(\A)=\pd\,\A$ by Corollary \ref{equa-ineq} (a). 

If $\pd_\A\,(\A)=\infty$, then by 
Corollary \ref{equa-ineq} (b) we get that $\coresdim_\B\,(\A)=\infty$, thus $\coresdim_\B\,(\mathrm{Mod}\,R)=\infty$, and since $\omega\subseteq \B$, also $\coresdim_\omega\,(\A)=\infty.$

If $\pd_\A\,(\A)<\infty$, then Corollary \ref{equa-ineq} (b) yields $\mathrm{Mod}\,R=\B^\vee$,
and we can apply Proposition \ref{Prop1} (b) since $\id_\A\,(\B)=0$.
So $\coresdim_\B\,(\mathrm{Mod}\,R)=\id_\A\,(\mathrm{Mod}\,R)=\pd\,\A.$ 
Further $\A\subseteq \omega^\vee$ by Lemma \ref{equal} (a). So, as in \cite[Corollary 2.3]{MS},  we infer from Proposition \ref{Prop1}   that 
$\coresdim_\B\,(\A)=\coresdim_\omega\,(\A)=\id_\omega \,(\A)=\id_\A\,\A=\pd\,\A.$
\

(b)   $\Rightarrow$: We have seen above  that $\A\subseteq\omega^\vee$ and  $\mathrm{Mod}\,R=\B^\vee$. Then the  result dual to Lemma  \ref{equal} (c) yields  that $\coresdim_\B\,(\mathrm{Mod}\,R)=\id_\omega\,(\mathrm{Mod}\,R)=\pd\,\omega.$ Hence we conclude  as in (a) that $\pd\,\omega=\id_\omega \,(\A)=\pd_\A\,(\A)<\infty.$
\

$\Leftarrow:$ Since $\omega\subseteq\A\subseteq\omega^\vee,$  we deduce from \ref{basic} (b) that $\pd\,\A
=\pd\,\omega<\infty.$  
\end{dem}

\bigskip

 \section{Applications to tilting theory}

We now apply our previous results to cotorsion pairs arising from tilting theory. 
\
 
\begin{defi}\cite{AnC} A module $T\in\mathrm{Mod}\,R$ is called tilting module provided
 \begin{itemize}
  \item[(T1)] $\pd\,T<\infty,$
\vspace{.2cm}
  \item[(T2)] $\Ext^i_R(T,T^{(I)})=0$ for any $i>0$ and all sets $I,$
\vspace{.2cm}
  \item[(T3)] $\coresdim_{\Add\,T}(R_R)<\infty.$
 \end{itemize}
\end{defi}

Every tilting module $T$ induces a complete hereditary cotorsion pair 
 $(\A,\B)$ in $\mathrm{Mod}\,R$ where $\B=T^\perp.$ It is  called the {\bf{tilting cotorsion pair}} induced by $T.$ 
 Tilting cotorsion pairs are characterized in \cite{AnC} by the property that the kernel is closed under coproducts and $\pd\,\A<\infty.$  We can now refine this characterization  as follows.
\

\begin{teo}\label{tiltingcp} Let $(\A,\B)$ be a hereditary complete cotorsion pair in $\mathrm{Mod}\,R$ with kernel  $\omega$. 
 \begin{itemize}
   \item[(a)] The following conditions are equivalent:
    \begin{itemize}
\vspace{.2cm}
     \item[(a1)] $(\A,\B)$ is a tilting cotorsion pair.
\vspace{.2cm}
     \item[(a2)] $\pd_\A\,(\A)<\infty$ and   $\omega$ is closed under coproducts.
\vspace{.2cm}
     \item[(a3)]    $\omega$ is closed under coproducts,  $\pd\,\omega<\infty$ and $\A\subseteq\omega^\vee$.
\vspace{.2cm}
    \end{itemize}
   \item[(b)] If $(\A,\B)$ is a tilting cotorsion pair induced by a tilting $R$-module $T,$ then
$$\pd\,T=\pd_\A\,(\A)=\coresdim_\B\,(\A)=\coresdim_{\Add\,T}\,(\A)=\coresdim_\B\,(\mathrm{Mod}\,R).$$
 \end{itemize} 
\end{teo}
\begin{dem} By \cite[Theorem 4.1]{AnC}, we know that $(\A,\B)$ is a tilting cotorsion pair if and only if $\pd\,\A<\infty,$  and    $\omega$ is closed under coproducts. So, the result follows from Theorem \ref{Thmequal} using that  $\omega=\Add\,T$, cf.~\cite[Lemma 2.4]{AnC}. 
\end{dem}

\

An important result of P. Eklof and J. Trlifaj \cite[Theorem 10]{ET} states   that  any set  $\X$ of $R$-modules gives rise to a complete cotorsion pair $({}^{\perp_1}(\X^{\perp_1}),\X^{\perp_1})$, known as the cotorsion pair {\bf{generated}} by $\X.$ As a consequence, every resolving subcategory  $\X$ of $\mathrm{mod}\,R$ consisting of modules of bounded projective dimension generates a tilting cotorsion pair, see \cite[Theorem 2.2]{AnHT}. Also this  result can be   reduced to the relative  projective dimension.
\

\begin{cor}\label{Corotiltingcp} Let $\mathcal{D}$ be a resolving subcategory of $\mathrm{mod}\,R$ such that $\pd_{\mathcal{D}}\,(\mathcal{D})<\infty$.
Then   the cotorsion pair $(\A, \B)$ generated by $\mathcal{D}$ is a  tilting cotorsion pair.
\end{cor}
\begin{dem} We know that $(\A, \B)$ is hereditary, complete, and $\B$ is closed under coproducts. Then also the kernel $\omega$ is closed under coproducts. Moreover, $\A\subseteq\varinjlim\mathcal{D}$ by \cite[Theorem 2.3]{AnT2}, so every $A\in\A$ has the form $A=\varinjlim D_i.$ Hence, for any $D\in\mathcal{D}$ and $j>n,$ we have
$$\Ext_R^j(D,A)=\Ext_R^j(D,\varinjlim D_i)\simeq \varinjlim\,\Ext_R^j(D,D_i)=0$$ since $D\in\mathrm{mod}\,R$ and $\pd_{\mathcal{D}}\,(\mathcal{D})=n.$ In other words, $\mathcal{D}\subseteq {}^{\perp_j}\A$ for any $j>n.$ Note that $\A$ consists  of all direct summands of $\mathcal{D}$-filtered modules, see \cite[3.2.3]{GT}. So,
it follows from  \cite[7.3.4]{EJ} that  $\A\subseteq {}^{\perp_j}\A$ for $j>n,$ proving that $\pd_\A\,(\A)\leq n.$  Now the statement follows immediately from \ref{tiltingcp}.
\end{dem}

\bigskip

\section{Finitistic dimensions}

Let us now consider a slightly more general situation. We assume that our cotorsion pair  $(\A,\B)$ satisfies  $\pd_\A\,(\A)<\infty$, but we drop the assumption that $\omega$ is closed under coproducts. It turns out that this is the setup needed for studying the finitistic dimensions.
\

We start with some general properties.

\begin{pro}\label{Proprincipal} Let $(\A,\B)$ be a hereditary complete cotorsion pair in $\mathrm{Mod}\,R$ with kernel $\omega.$ Assume $\pd_\A\,(\A)<\infty$. The following statements hold true.
 \begin{itemize}
\vspace{.2cm}
  \item[(a)] $\A^\wedge=\{M\in\mathrm{Mod}\,R\,:\,\pd_\B\,(M)<\infty\}=\Q.$
\vspace{.2cm}
  \item[(b)] $\resdim_\A\,(M)=\pd_\B\,(M)$ for any $M\in\Q.$
\vspace{.2cm}
  \item[(c)] $\pd\,M\leq\pd_\A\,(\A)+\resdim_\A\,(M)$ for any $M\in\Q.$
  \vspace{.2cm}
  \item[(d)] $\Q\cap\B=\omega^\wedge.$
 \end{itemize}
\end{pro}
\begin{dem} Since $\pd_\B\,(\A)=0$ and $\A={}^{\perp_1}\B,$ we get from the statement   dual to  \ref{Prop1} (b) that $\pd_\B\,(M)=\resdim_\A\,(M)$ for any $M\in\A^\wedge.$ This yields the inclusion $\A^\wedge\subseteq\{M\in\mathrm{Mod}\,R\,:\,\pd_\B\,(M)<\infty\}.$ The reverse inclusion follows from the statement  dual to \ref{Prop1} (c). On the other hand, since $\Q_0\subseteq\A,$ we have that $\resdim_\A\,(M)\leq\resdim_{\Q_0}\,(M)=\pd\,M$ for any $M\in\mathrm{Mod}\,R,$ implying $\Q\subseteq\A^\wedge.$ Finally,  $\pd\,\A<\infty $ by Corollary \ref{equa-ineq}. Then we infer from Lemma \ref{basic} (c) that $\A^\wedge\subseteq\Q,$ which completes the proof of (a) and (b).
\\
Dualizing  \ref{Prop1} (a), we obtain $\pd_\A\,(M)\leq\pd_\A\,(\A)+\resdim_\A\,(M)$ for any $M\in\A^\wedge=\Q.$ Note moreover that $\Q=\Q_0^\wedge$ and  $\pd\,M=\resdim_{\Q_0}\,(M)$. Thus
$\pd_\A\,(M)=\pd\,M$ by the dual of \ref{Prop1} (b).
This proves  (c). 
\\
Finally, we prove (d). We have that $\omega^\wedge\subseteq\Q\cap\B$ since $\Q$ and $\B$ are closed under cokernels of monomorphisms, and $\pd\,\omega<\infty$ by Theorem \ref{Thmequal}. For the reverse inclusion,  let $X\in\Q\cap\B$ and $r=\pd_\B\,(X)\leq\pd\,X.$ Since $(\A,\B)$ is complete and $\A\subseteq\Q$, there is an exact sequence $0\to K_r\to W_r\stackrel{f_r}{\to}W_{r-1}\to\cdots\to W_1\stackrel{f_1}{\to}W_0\stackrel{f_0}{\to}X\to 0$ such that $W_i\in\omega$ and $K_i:=\Ker(f_i)\in\Q\cap\B$ for each $i.$ By assumption, $\Ext_R^{r+1}(X,-)|_\B=0$. Since $W_i\in{}^\perp K_{r+1},$ we infer by dimension shifting that $0=\Ext_R^{r+1}(X,K_r)\simeq\Ext_R^1(K_{r-1},K_r).$ Hence $K_{r-1}\in\omega,$ and $X\in\omega^\wedge$. 
\end{dem}
\

\

We now introduce two numbers $\alpha=\resdim_\A\,(\Q)$ and $\beta=\resdim_\A\,(\Q^{<\infty})$  which allow to obtain interesting estimates for the finitistic dimensions.
\

\begin{lem}\label{specialnum} Let $(\A,\B)$ be a hereditary complete cotorsion pair in $\mathrm{Mod}\,R$ with kernel  $\omega$. Assume $\pd_\A\,(\A)<\infty$, and set $\alpha=\resdim_\A\,(\Q)$. The following statements hold true.
\begin{enumerate}
\item[(a)] $\alpha=\resdim_\A\,(\omega^\wedge)=\resdim_\omega\,(\omega^\wedge)=\pd_\omega\,(\Q)\leq\pd\,\omega^\wedge.$
\vspace{.2cm}
\item[(b)] $\pd\,\omega^\wedge=\Findim(R_R).$
\end{enumerate}
\end{lem}
\begin{dem}
(a) If   $M\in\Q$, then $\resdim_\A\,(M)=\pd_\B(M)$ by Proposition \ref{Proprincipal} (b).  Further, taking an exact sequence $0\to M\to B\to A\to 0$ where $B\in\B$ and $A\in\A,$ we have $B\in\Q\cap\B=\omega^\wedge.$ Thus we infer from  Lemma \ref{basic} (b) that
  $\pd_\B(M)\leq\pd_\B\,(B)\leq\pd_\B\,(\omega^\wedge)$.
  This proves $\alpha\le \pd_\B\,(\omega^\wedge)$.
  On the other hand, $\resdim_\A\,(\omega^\wedge)\le \alpha$ since $\omega^\wedge\subseteq\Q$.
  Hence $\alpha=\resdim_\A\,(\omega^\wedge)=\pd_\B\,(\omega^\wedge)\le \pd\,(\omega^\wedge)$ by Proposition \ref{Proprincipal}.\\
  Finally, by Lemma \ref{equal} (b) and (c) we have $\alpha=\resdim_\omega\,(\omega^\wedge)=\pd_\omega\,(\Q)$.
\

(b)  We know from Proposition \ref{Proprincipal}   that    $\omega^\wedge=\Q\cap\B\subseteq \Q$, so $\pd\,\omega^\wedge\leq\Findim(R_R).$ 
Now let $M\in\Q$. As above, there is   an exact sequence $0\to M\to B\to A\to 0$ where $B\in\Q\cap\B=\omega^\wedge$ and $A\in\A.$ Then we know by Lemma \ref{basic} (b) that
$\pd\,M\leq\mathrm{max}\{\pd\,B,\,\pd\,A-1\}$.
But $\pd\,\A=\pd\,\omega$ by Theorem \ref{Thmequal}, so $\pd\,A\leq\pd\,\omega^\wedge$.
Thus $\pd\,M\leq\pd\,\omega^\wedge$, and we conclude  $\Findim(R_R)=\pd\,\omega^\wedge.$ 
\end{dem}
\

\begin{teo}\label{Thmprincipal} Let $(\A,\B)$ be a hereditary complete cotorsion pair in $\mathrm{Mod}\,R$. Assume  $\pd_\A\,(\A)<\infty,$ and set  $\alpha=\resdim_\A\,(\Q)$, and $\beta=\resdim_\A\,(\Q^{<\infty})$. The following statements hold true.
\begin{itemize}
 \item[(a)] $\alpha\leq\Findim(R_R)\leq\pd_\A\,(\A)+\alpha.$
\vspace{.2cm}
 \item[(b)] $\beta\leq\findim(R_R)\leq\pd_\A\,(\A)+\beta.$ 
\end{itemize}
\end{teo}
\begin{dem} (a) $\alpha\leq\pd\,\omega^\wedge=\Findim(R_R)$ by Lemma \ref{specialnum}, and $\Findim(R_R)\leq\pd_\A\,(\A)+\alpha$ by Proposition \ref{Proprincipal} (c).
\

(b) By Proposition \ref{Proprincipal} (b), we have $\beta=\pd_\B\,(\Q^{<\infty})\leq\pd\,\Q^{<\infty}=\findim(R_R).$ Furthermore,  $\findim(R_R)\leq\pd_\A\,(\A)+\beta$ by Proposition \ref{Proprincipal} (c).
\end{dem}
\

\

We now obtain criteria for finiteness of the finitistic dimensions.

\begin{cor}\label{equivFindim} Let $R$ be a ring.
\begin{enumerate}
\item[(I)] The following conditions are equivalent.
\vspace{.2cm}
 \begin{itemize}
  \item[(a)] $\Findim(R_R)<\infty.$
\vspace{.2cm}
  \item[(b)] There is a hereditary complete cotorsion pair $(\A, \B)$ in $\mathrm{Mod}\,R$ such that $\pd_\A\,(\A)<\infty$ and $\alpha=\resdim_\A\,(\Q)<\infty.$
\vspace{.2cm}
  \item[(c)] Every hereditary complete cotorsion pair $(\A, \B)$ in $\mathrm{Mod}\,R$ such that $\pd_\A\,(\A)<\infty$ satisfies $\alpha=\resdim_\A\,(\Q)<\infty.$
 \end{itemize}
 
 \medskip
 
 \item[(II)]  The following conditions are equivalent.
\vspace{.2cm}
 \begin{itemize}
  \item[(a)] $\findim(R_R)<\infty.$
\vspace{.2cm}
  \item[(b)] There is a hereditary complete cotorsion pair $(\A, \B)$ in $\mathrm{Mod}\,R$ such that $\pd_\A\,(\A)<\infty$ and $\beta=\resdim_\A\,(\Q^{<\infty})<\infty.$
\vspace{.2cm}
  \item[(c)] Every hereditary complete cotorsion pair $(\A, \B)$ in $\mathrm{Mod}\,R$ such that $\pd_\A\,(\A)<\infty$ satisfies $\beta=\resdim_\A\,(\Q^{<\infty})<\infty.$
 \end{itemize}
 \end{enumerate}
\end{cor}
\begin{dem} 
(I)
The implications (b)$\Rightarrow$(a) and (a)$\Rightarrow$(c) follow  directly from Theorem \ref{Thmprincipal}. To prove (c)$\Rightarrow$(b), we choose the cotorsion pair $(\Q_n,\Q_n^\perp)$ for some natural number $n.$ It is well known that $(\Q_n,\Q_n^\perp)$ is hereditary   and complete, see  \cite[Theorem 7.4.6]{EJ}. Furthermore $\pd_{\Q_n}\,(\Q_n)\leq n$ by Corollary \ref{equa-ineq}, so $\alpha$ is finite  by the assumption in (c).
\

(II) It is proven similarly.
\end{dem}

\medskip

Corollary \ref{equivFindim} (II)   characterizes the  validity of the Second Finitistic Dimension Conjecture. The First Finitistic Dimension Conjecture stated that $\Findim(R_R)=\findim(R_R)$. It is well known  by now   that this conjecture fails in general. The first example of this phenomenon was given by Huisgen-Zimmermann
in 1992.
Later, Smal\o\ even showed that the difference between the big and the little finitistic dimension can be arbitrarily large, see \cite{Sm}.
However, in many interesting cases the two numbers do coincide. The following result
characterizes   this situation.
\

\begin{teo}\label{Findim-findim} Let $R$ be a right noetherian ring with $\findim(R_R)<\infty$. Let $(\A, \B)$ be the cotorsion pair generated by $\Q^{<\infty},$ and let $\omega$ be its kernel. The following statements hold true.
 \begin{itemize}
  \item[(a)] $\Findim(R_R)-\findim(R_R)\leq\alpha=\resdim_\A\,(\Q).$
\vspace{.2cm}
  \item[(b)] $\Findim(R_R)=\findim(R_R)$ if and only if $\pd\,\omega^\wedge=\pd\,\omega.$
 \end{itemize}
\end{teo}
\begin{dem} We know from \cite[Theorem 2.6]{AnT1} and Theorem \ref{tiltingcp} (b) that $(\A, \B)$ is a tilting cotorsion pair with $\pd_\A\,(\A)=\findim(R_R).$ So, statement (a)   follows from Theorem \ref{Thmprincipal}. On the other hand, by Proposition \ref{Thmequal} and Lemma \ref{specialnum} (b) we further have $\findim(R_R)=\pd\,\omega$, and $\Findim(R_R)=\pd\,\omega^\wedge$, proving statement (b).
\end{dem}
\

\begin{rk}{\rm Let  $(\A, \B)$ be a cotorsion pair in $\mathrm{Mod}\,R,$ and set $\alpha=\resdim_\A\,(\Q)$.

 \medskip

$(1)$ Note that $\alpha<\infty$ does not imply that $\pd_\A\,(\A)<\infty.$ For example, if we  take the trivial cotorsion pair given by $\A=\mathrm{Mod}\,R$ and $\B=\I_0,$  we always have $\alpha=0$ since $\Q\subseteq \A$, but $\pd_\A\,(\A)$ is the global dimension of $R$.
\

\medskip

 $(2)$ Let $(\A, \B)$ be hereditary and complete with $\pd_\A\,(\A)<\infty$. Then   $\alpha=0$ if and only if $\A=\Q$, and in this case  $\Findim(R_R)=\pd_\A\,(\A)<\infty.$
 
In particular, if $(\A, \B)$ is the cotorsion pair generated by $\Q^{<\infty}$, then $\alpha=0$ implies that  $\Findim(R_R)=\findim(R_R)=\pd_\A\,(\A)<\infty.$ 

So, our results generalize
 \cite[Proposition 3.1 and Theorem 3.2]{AnT1}.
 
 \medskip
 
 $(3)$ The inequality in Theorem \ref{Findim-findim} (a) can be strict. For example, if $R$ is the finite dimensional algebra
from \cite{IST}, then $\Findim(R_R)=\findim(R_R)=1,$ but  the cotorsion pair   $(\A, \B)$  generated by $\Q^{<\infty}$ satisfies $\alpha=1.$ In fact, it is shown in \cite{AnT2} that  $\A\subsetneq\Q$, thus    $0<\alpha\leq\Findim(R_R)=1$ by Theorem \ref{Thmprincipal}.

\medskip

$(4)$ Let $R$ be a ring with   $\findim(R_R)<\infty.$ If $(\A, \B)$ is the cotorsion pair generated by $\Q^{<\infty},$ then  $\Q^{<\infty}\subseteq\A,$ and so $\beta=\resdim_\A\,(\Q^{<\infty})=0.$
}
\end{rk}

\

We close the paper with a result  concerning arbitrary tilting modules.

\

\begin{cor}\label{pd+id} Let $T$ be a tilting module in $\mathrm{Mod}\,R.$ If $R$ is right noetherian then $$\Findim(R_R)\leq\pd\,T+\id\,T.$$
\end{cor}
\begin{dem} 
We apply the first formula in  Theorem \ref{Thmprincipal} to the tilting cotorsion pair $(\A,\B)$ induced by $T$. Recall that its kernel $\omega$ coincides with $\Add T$, and $\pd_\A\,(\A)=\pd\,T$ by Theorem \ref{tiltingcp}. So,  we obtain $\Findim(R_R)\leq\pd\, T+\alpha,$ where $\alpha=\pd_{\Add\,T}\,(\Q)=\id_\Q\,(\Add\,T)\leq\id\,(\Add\,T)$ by Lemma \ref{specialnum} (a). Now, if $R$ is right noetherian, then
$\id\,(\Add\,T)=\id\,T$.
\end{dem}
\

\begin{rk}\label{applications} {\rm Let $R$ be right noetherian.

\medskip

$(1)$  Note that there is a 
tilting module $T$   in $\mathrm{Mod}\,R$ with $\id\,T<\infty$  if and only if $\id\,R_R<\infty$. Indeed,   the only-if-part holds true by condition (T3) in the definition of a tilting module, and  for the if-part  one chooses $T=R_R$.
\

Choosing $T=R_R$, Corollary \ref{pd+id} gives $\Findim(R_R)\leq\id\,(R_R).$ This was proved by H. Bass in \cite{B}.

\medskip

 $(2)$ The inequality in \ref{pd+id} is sharp. 
 To see this, let $R$ be an Iwanaga-Gorenstein ring with $\id\,(R_R)=n.$ Consider the minimal injective coresolution $0\to R_R\to I_0\to I_1\to\cdots\to I_n\to 0.$ By \cite[Section 3]{AnHT} the injective module $T:=\oplus_{i=0}^n\,I_i$ is tilting, and $\Findim(R_R)=n=\pd\,T=\pd\,T+\id\,T$. 
We could also choose the tilting   module $T=R_R$. In this case   $\Findim(R_R)=\id\,T=\pd\,T+\id\,T.$

\medskip

 $(3)$ In case $R$ is an artin algebra, $T={}_RT$ is a classical tilting and cotilting left $R$-module and $A:=\mathrm{End}(T),$ it was proven in \cite{W} that $\findim({}_AA)\leq\pd\,T+\id\,T.$

\medskip

 $(4)$ If $R$ is a properly stratified algebra having a simple preserving duality, and such that every classical tilting right $R$-module is also cotilting (see \cite{MO}), then $\Findim(R_R)=\findim(R_R).$ Indeed, by \cite[Theorem 1]{MO}, we have that $\findim(R_R)=2\,\pd\,T$ where $T$ is the so called characteristic tilting module associated to the properly stratified algebra $R.$ Moreover, as a consequence of the simple preserving duality, it can be seen that $\pd\,T=\id\,T.$ Therefore, from \ref{pd+id} it follows that $2\,\pd\,T=\findim(R_R)\leq\Findim(R_R)\leq 2\,\pd\,T,$ proving our assertion.
} 
\end{rk}

\bigskip

{\bf Acknowledgements.} We thank Steffen K\"onig for drawing our attention on the paper \cite{MO}.
This research was carried out while the second named  author was a Visiting Professor of the Facolt\`a di Scienze dell'Universit\`a dell'Insubria at Varese; he enjoyed his visit very much, and he would like to thank  for the hospitality. The first named author was partially supported by the DGI and the European Regional
Development Fund, jointly, through Project
 MTM2005--00934, and by the Comissionat per Universitats i Recerca
of the Generalitat de Ca\-ta\-lunya, Project 2005SGR00206, and by Progetto di Ateneo CPDA071244 of the University of Padova. The second named author was partially supported by Project PAPIIT-Universidad Nacional Aut\'onoma de M\'exico IN101607.

\footnotesize

\vskip3mm \noindent Lidia Angeleri H\"ugel\\ Dipartimento di Informatica
e Comunicazione\\
Universita degli Studi dell'Insubria.\\
Via Mazzini 5, I-21100. Varese, ITALY.\\ {\tt lidia.angeleri@uninsubria.it}

\vskip3mm \noindent Octavio Mendoza Hern\'andez\\ Instituto de Matem\'aticas\\ Universidad Nacional Aut\'onoma de M\'exico.\\ 
Circuito Exterior, Ciudad Universitaria\\
C.P. 04510, M\'exico, D.F. MEXICO.\\ {\tt omendoza@matem.unam.mx}


\bigskip


\begin{thebibliography}{20}
\bibitem{AnC} L. Angeleri H\"ugel, F.U. Coelho. Infinitely generated tilting
modules of finite projective dimension. {\it Forum Math.}13(2001)239-250.
\bibitem{AnT1} L. Angeleri H\"ugel, J. Trlifaj. Tilting Theory and the Finitistic
Dimension Conjectures. {\it Trans. of the AMS.} Vol. 354(2002)4345-4358.
\bibitem{AnT2} L. Angeleri H\"ugel, J. Trlifaj. Direct limits of modules of 
finite projective dimension. {\it Rings, Modules, Algebras, and Abelian Groups, LNPAM}
236, M.Dekker, New York (2004), 27-44.
\bibitem{AnHT} L. Angeleri H\"ugel, D. Herbera, J. Trlifaj. Tilting modules and Gorenstein rings.
{\it Forum Math.} 18(2006)211-229.
\bibitem{AuB} M. Auslander, R.O. Buchweitz.  The Homological Theory of
  Maximal Cohen-Macaulay Approximations. {\it Societe Mathematique de
  France}. Memoire 38 (1989) 5-37.
\bibitem{B} H. Bass. Injective dimension in noetherian rings. {\it Trans. Amer. Math. Soc.}102(1962)18-29.
\bibitem{EJ} E. Enochs, O.M.G. Jenda. Relative Homological Algebra. 
{\it De Gruyter, Berlin.}(2001). 
\bibitem{GT} R. Gobel, J. Trlifaj. Approximations and Endomorphism Algebras
of Modules. {\it Editorial Walter de Gruyter.} (2006).
\bibitem{ET} P.Eklof, J. Trlifaj. How to make Ext vanish. {\it Bull. London Math. Soc.}
33(2001)41-51.
\bibitem{IST} K. Igusa, S.O. Smalo, G. Todorov. Finite projectivity and contravariant finiteness.
{\it Proc. Amer. Math. Soc.} 109(1990)937-941.
\bibitem{MS} O. Mendoza, C. Saenz. Tilting Categories with applications to 
Stratifying Systems. {\it Journal of algebra.} 302(2006)419-449. 
\bibitem{MO} V. Mazorchuk, S. Ovsienko. Finitistic dimension of properly stratified algebras. {\it Adv. in Math.} 186 (2004) 251-265.
\bibitem{S} L. Salce. Cotorsion theories for abelian groups. {\it{Symposia Math.}}XXIII(1979)
11-32.
\bibitem{Sm}{S.O.Smal\o}, Homological difference between finite and
infinite dimensional representations of algebras, Trends in Math.,
Birkh\" auser, Basel 2000, 81-93.
\bibitem{W} J. Wei. Finitistic Dimension and Restricted Flat Dimension.
{\it Journal of algebra} Volume 320, Issue 1 (2008) Pages 116-127.
\end{thebibliography}
\end{document}